\documentclass{amsart}

\usepackage{color}
\usepackage[colorlinks=true, citecolor=blue, linkcolor=blue]{hyperref}

\usepackage{amsmath}

\usepackage{amssymb}

\usepackage[round]{natbib}

\usepackage{xparse}
\NewDocumentCommand{\ceil}{s O{} m}{%
  \IfBooleanTF{#1} 
    {\left\lceil#3\right\rceil} 
    {#2\lceil#3#2\rceil} 
}
\NewDocumentCommand{\floor}{s O{} m}{%
  \IfBooleanTF{#1} 
    {\left\lfloor#3\right\rfloor}
    {#2\lfloor#3#2\rfloor}
}

\definecolor{c20}{rgb}{0.,0.7,0.}
\definecolor{c30}{rgb}{0.,0.,1.}
\definecolor{c40}{rgb}{1,0.1,0.7}
\definecolor{c50}{rgb}{1,0,0}
\definecolor{c60}{rgb}{1,0.9,0.1}
\definecolor{c90}{cmyk}{0, 1, 0.1, 0.5}

\def\BH#1{\textcolor{c20}{#1}}
\def\BH#1{#1}

\def\KD#1{\textcolor{c90}{#1}}
\def\KD#1{#1}

\def\He#1{\textcolor{c50}{#1}}
\def\He#1{#1}
\def\Eb#1{\textcolor{c50}{#1}}
\def\Eb#1{#1}
\def\EEH#1{\textcolor{c50}{#1}}
\def\EEH#1{#1}
\def\Ec#1{\textcolor{c50}{#1}}
\def\Ec#1{#1}

\newtheorem{cond}{Condition}
\newtheorem{lemma}{Lemma}
\newtheorem{example}{Example}
\newtheorem{theorem}{Theorem}
\newtheorem{remark}{Remark}

\newcommand{\ve}{\varepsilon}

\newcommand{\abs}[1]{\left\lvert #1 \right\rvert}

\newcommand{\E}[1]{\mathbb{E}\left\{#1\right\}}

\newcommand{\pk}[1]{\mathbb{P} \left\{ #1 \right\} }

\newcommand{\R}{\mathbb{R}}
\newcommand{\Z}{\mathbb{Z}}

\newcommand{\N}{\mathbb{N}}
\newcommand{\inr}{\in \R}
\newcommand{\inn}{\in \N}
\newcommand{\ldot}{,\ldots,}

\newcommand{\limit}[1]{\lim_{#1 \to   \infty}}

\newcommand{\BQN}{\begin{eqnarray}}
\newcommand{\EQN}{\end{eqnarray}}
\newcommand{\BQNY}{\begin{eqnarray*}}
\newcommand{\EQNY}{\end{eqnarray*}}

\renewcommand{\P}{\mathbb P}
\newcommand{\BS}{\begin{sat}}
\newcommand{\ES}{\end{sat}}
\newcommand{\BT}{\begin{theorem}}
\newcommand{\ET}{\end{theorem}}
\newcommand{\BK}{\begin{korr}}
\newcommand{\EK}{\end{korr}}
\newcommand{\EQD}{\stackrel{d}{=}}

\newcommand{\BD}{\begin{de}}
\newcommand{\ED}{\end{de}}

\newcommand{\BIT}{\begin{itemize}}
\newcommand{\EIT}{\end{itemize}}
\newcommand{\BDI}{\begin{description}}
\newcommand{\EDI}{\end{description}}

\newcommand{\BRM}{\begin{remark}}
\newcommand{\ERM}{\end{remark}}

\newcommand{\BEL}{\begin{lemma}}
\newcommand{\EEL}{\end{lemma}}

\newcommand{\nelem}[1]{{Lemma \ref{#1}}}

\newcommand{\netheo}[1]{{Theorem \ref{#1}}}

\newcommand{\prooftheo}[1]{ \textsc{\bf Proof of Theorem} \ref{#1}:}

\newcommand{\prooflem}[1]{\textsc{\bf Proof of Lemma} \ref{#1}:}

\newcommand{\COM}[1]{}

\newcommand{\QED}{\hfill $\Box$}

\def\xiW{\xi_W}

\def\IF{\infty}

\date{}

\def\Var{\text{Var}}

\def\HWD{\mathcal{H}_W^\delta}

\def\md{M^\delta}

\begin{document}

\title[Generalized Pickands constants]{Generalized Pickands constants and stationary max-stable processes}

\author{Krzysztof D\c{e}bicki}
\address{Krzysztof D\c{e}bicki, Mathematical Institute, University of Wroc\l aw, pl. Grunwaldzki 2/4, 50-384 Wroc\l aw, Poland}
\email{Krzysztof.Debicki@math.uni.wroc.pl}
\author{Sebastian Engelke}
\address{Sebastian Engelke, Ecole Polytechnique F\'ed\'erale de Lausanne,
EPFL-FSB-MATHAA-STAT, Station 8, 1015 Lausanne, Switzerland
}
\email{sebastian.engelke@epfl.ch}
\author{Enkelejd Hashorva}
\address{Enkelejd Hashorva, University of Lausanne, B\^{a}timent Extranef, UNIL-Dorigny, 1015 Lausanne, Switzerland}
\email{Enkelejd.Hashorva@unil.ch}

\begin{abstract}
Pickands constants play a crucial role in the
asymptotic theory of Gaussian processes. They are
commonly defined as the limits of a sequence of expectations involving fractional Brownian motions
and, as such, their exact value is often unknown. Recently,
\cite{DiekerY} derived a novel representation of Pickands constant
as a simple expected value that does not involve a limit operation.
In this paper we show that
the notion of Pickands constants and their
corresponding Dieker--Yakir representations can be
extended to a large class of stochastic processes,
including general Gaussian and L\'evy processes.
We furthermore provide a link to spatial extreme value
theory and show that Pickands-type  constants coincide
with certain constants arising in the study of max-stable
processes with mixed moving maxima representations.
\end{abstract}

\maketitle

\section{Introduction}
\KD{
The asymptotic behavior of the tail probabilities of the supremum of a Gaussian process $\{X(t),t\in[0,T]\}$, $T>0$, with continuous sample paths
is well understood for a wide class of correlation structures
of $X$.
Its general form, that is valid for both
the classical
Pickands' theorem for the centered stationary case and the result by Piterbarg for the non-stationary case, shows that for any
$\delta \ge 0$ (set $\delta \Z= \R$ if $\delta=0$)
\BQN
\pk{\sup_{t\in \delta \Z \cap [0,T]} X(t)> u} \sim  \HWD  C u^a e^{- u^2/b}, \quad u\to \IF,
\EQN
holds under some mild regularity conditions on the correlation and the variance function of $X$
\citep{PickandsB,Berman82,Berman92,Pit96,Pit20,Tabis}.
Here all the positive constants $a, b,C$ are explicitly known,}
whereas the constant $\HWD$,
which is referred to as Pickands constant, is given by the following limit
 \BQN\label{pickands_const}
 \HWD
 = \limit{T} T^{-1} \E{ \sup_{t\in \delta \Z  \cap [0,T] } e^{ W(t)}} \in(0,\IF), \quad W(t)=\sqrt{2} B_{\alpha}(t)- \abs{t}^\alpha,
\EQN
where $\{B_\alpha(t),t\ge 0\}$ is a centered  fractional Brownian motion
with Hurst index $\alpha/2 \in (0,1]$, that is, a mean zero Gaussian processes  with continuous sample paths and
covariance function
$$ \text{Cov}\{B_\alpha(s), B_\alpha(t)\}= \frac{1}{2} \Bigl( \lvert t\rvert^\alpha + \lvert s\rvert^\alpha- \lvert t-s\rvert^\alpha\Bigr),  s,t\ge 0.$$
The only known values of $\HWD$ are for $\delta=0$ if $\alpha=1,2$.
Numerous papers have considered the calculation of Pickands constants, \KD{with particular focus on the case}
$\delta=0$; see for instance \citet{shao1996bounds, MR1708208, demiro2003simulation, Krzys2006Pickands, debicki2008note, Harper1, Harper2}. \\
Recently, the seminal contribution \cite{DiekerY} derived the alternative representation for $\HWD$
\BQN\label{repDY}
\HWD=\E{ \frac{\md}{ S^\eta }}, \quad \forall \delta=\eta > 0, \text{  or } \delta=0, \eta \ge 0,
\EQN
{\KD{where}
\BQN\label{md}
 \md = \sup_{t\in \delta \Z  } e^{W (t)}, \quad S^\eta =\eta \sum_{t\in \eta \Z} e^{W(t)},  \quad S^0 = \int_{\mathbb R} e^{W(t)}\, dt.
 \EQN
\BH{The principal advantage of
\KD{{\it Dieker--Yakir representation}}
\eqref{repDY}} is that it is given as an expectation rather than as a limit, which is particularly useful for Monte Carlo simulations of $\HWD$. \\
\BH{Pickands constants traditionally also appear in Gumbel limit theorems, see e.g., \cite{Berman92, Pit2004}. Such limit theorems
are recently formulated for max-stable processes and provide a first link of classical
Gaussian tail asymptotics to spatial extreme value theory. Specifically,
\cite{DM} showed that \cite[see also][]{stoev2010max, MIKYuw}}
\BQN \label{remarkG}
\limit{T} \pk{ \sup_{t\in \delta \Z \cap [0,T]} \xiW (t) \le x + \ln T}
= \exp\left\{-  \HWD  \exp(-x)\right\}, \quad x\inr,
\EQN
where the so-called Brown--Resnick process $\xiW$ is defined as
\BQN\label{eA}
 \xiW (t)= \max_{i\ge 1}( P_i+ W_i(t)), \quad t\inr.
\EQN Here $\Pi=\sum_{i=1}^\IF \ve_{P_i}$ is a Poisson point process
with intensity $e^{-x} dx$, and $W_i,i\ge 1$, are independent  copies of $W$, also independent of $\Pi$. We denote by $\ve_x$ the unit Dirac measure at $x\in\R$.
The Brown--Resnick process $\xiW$ is both max-stable and stationary \citep{kab2009, kab2009a, kab2011, MolchanovSPA, MolchanovBE}.
The stationarity means that the processes $\{ \xiW (t),t\inr\}$ and $\{ \xiW (t+h),t\inr\}$ have the same distribution for any $h\inr$. Moreover, the process $\xiW$  arises naturally as the limit of suitably normalized pointwise maxima
of independent copies of stationary Gaussian processes \citep[Theorem 17]{kab2009}. This makes this class of
processes a widely-used model in the risk assessment of spatial extreme events.
The result in \eqref{remarkG} states that $\HWD$ coincides with the so-called \emph{extremal index}
of the stationary, max-stable process $\xi_W$, a quantity that summarizes the
temporal extremal dependence \citep[c.f.,][]{lea1983}.\\
\COM{ A remarkable finding of \cite{DM} (see also  \cite{stoev2010max}, \cite{MIKYuw} for the Fr\'{e}chet  case) is the following Gumbel limit result
\BQN \label{remarkG}
\limit{T} \pk{ \sup_{t\in \delta \Z \cap [0,T]} \xiW (t) \le x + \ln T}
= \exp\left(-  \HWD  \exp(-x)\right), \quad x\inr.
\EQN}
Another interesting representation of Pickands constant for $\delta >0$ in the case of fractional Brownian motion in \cite{DiekerY} is
\BQN\label{imm2}
\HWD&=& \frac{1} \delta \pk{ \sup_{t\in \delta \Z } W(t)=0} =:C_W^\delta.
\EQN
\BH{Surprisingly, the constant $C_W^\delta$ appears} in the moving maxima representation of $\xiW$ restricted on $\delta \Z$; see Theorem 8 and Remark 9 in \cite{KabExt}. In the aforementioned contribution, the constant $C_W^\delta$ has already been evaluated numerically for different values of $\delta$ in order
to simulate samples from the max-stable process $\xiW$. This underlines the
connection between spatial extremes and classical asymptotic theory of Gaussian processes.

The objective of this paper is \KD{twofold.}
On the one hand, we consider generalized Pickands
constants $\mathcal{H}_W^\delta$ in \eqref{pickands_const}, where $W$ is replaced by
more general stochastic processes than fractional Brownian motions, which are not necessarily Gaussian.
We are then interested in finding conditions for the existence and positiveness of the limit in \eqref{pickands_const},
and in deriving equivalent representations \KD{of these constants}.
More precisely, we show that for $W$ chosen such that $\xiW$ is max-stable and stationary,
generalized Pickands constants can be defined in $(0,\infty)$,
and, most notably, that they admit a Dieker--Yakir type representation \eqref{repDY}
under certain conditions.\\
\Eb{On the other hand, we explore the connection between mixed moving maxima processes and generalized  Pickands constants that is suggested by equation \eqref{imm2}. Our findings are beneficial for both the theory of extremes of max-stable stationary processes, and the asymptotic theory of random processes. In particular, we show that $\HWD = C_W^\delta$, which holds
not only for $\delta>0$ but also in the classical case $\delta =0$. This shows that calculation of the classical Pickands constant is related to the simulation of the corresponding max-stable processes discussed above.
}

The paper is organized as follows. In Section \ref{sec_gen_pick}
we introduce generalized Pickands constants $\mathcal H_W^\delta$
and give conditions
under which they admit a Dieker--Yakir type representation.
Examples
for the process $W$ will be general Gaussian processes with stationary increments
and L\'evy processes. The connection of the constants $\mathcal H_W^\delta$
to mixed moving maxima processes is investigated
in Section \ref{sec_m3}. This link will provide a simple
proof of the positiveness of generalized Pickands constants.
All proofs are given in Section \ref{proofs}.
\KD{The Appendix} comprises some facts on discrete mixed moving maxima
representations which are needed in Section \ref{sec_m3}.

\section{Generalized Pickands constants}\label{sec_gen_pick}

Let $\{B(t), t\inr\}$ be a stochastic process
on the space $D$ of c\`adl\`ag functions $f:\R\to\R$ with $B(0) = 0$
and finite $\E{e^{B(t)}}< \infty$, for all $t\in\R$.
We introduce the drifted process
\BQN \label{Ws}
 W(t)= B(t) - \ln   \E{e^{B(t)}}, \quad t\inr
\EQN
and note that it satisfies $\E{e^{W(t)}}=1$.
We can
therefore define the corresponding max-stable process $ \xiW $ by the construction
\eqref{eA} which has standard Gumbel margins
Throughout, we will assume that $W$ is chosen such that
the process $ \xiW $ is stationary and has c\`adl\`ag sample paths; see Proposition 6 in \cite{kab2009} for a general stationarity criterion.

In this section we introduce the \emph{generalized Pickands constant}
of the process $W$ on the grid $\delta \Z$ for $\delta\geq 0$
as
\BQN\label{gen_pick}
\mathcal{H}_W ^\delta = \limit{T} \frac1T \E{ \sup_{t\in \delta \Z \cap [0,T] } e^{ W (t)}}.
\EQN
The existence of the expected value in \eqref{gen_pick} when $\delta=0$ is equivalent to the assumption that $ \xiW $ has c\`adl\`ag sample paths \citep{kabDombry}.
However, the existence and finiteness of the limit as $T\to\infty$ is not obvious. In the sequel, we investigate:
\begin{itemize}
  \item[a)]
    the existence of the constant $\mathcal{H}_W ^\delta$,
  \item[b)]
    its finiteness and positivity,
  \item[c)]
    equivalent representations that can for instance be used
    for efficient approximations.
\end{itemize}
In Section \ref{sec2.1} we discuss question a) in a general setting.
For question b) and c) we will concentrate on two important examples
for $W$ such that the above assumptions are satisfied.
In Section \ref{sec2.2} we consider the general Gaussian case, where
\begin{itemize}
  \item[$\diamond$]
    $B$ is a sample continuous centered
    Gaussian process
    with stationary increments and
    variance function $\sigma^2(t), t\inr$.  With
    $$W(t)= B(t) - \sigma^2(t)/2, \quad t\inr,$$
    the process $ \xiW $
    is max-stable and stationary. Its law depends only  on
    the variogram $\gamma(t)= \Var(B(t)- B(0))$ and
    we can therefore assume without loss of generality that $W(0)=0$;
    see \cite{kab2009, kab2011} for details.
\end{itemize}
The generalized Pickands constant can also be defined for non-Gaussian
processes. In Section \ref{sec2.3} we investigate the case where
\begin{itemize}
  \item[$\diamond$]
    $\{B^+(t), t\geq 0\}$ is a 
    L\'evy process such that
    $\Phi(\theta)=\ln \E{ e^{ \theta B^+(1)}}$
    is finite for $\theta =1$ and set
    $$W(t) = B^+(t) - \Phi(1)t, \quad t\geq 0.$$
    If $\{W(t), t\leq 0\}$ is defined as an
    exponentially transformed version of the corresponding $\{W(t), t\geq 0\}$, then
    $ \xiW $ can be shown to be stationary and max-stable; see \cite{sto2008,eng2014d} for details.
\end{itemize}
Clearly, these are not the only examples. For instance, a slight generalization is to introduce an independent mixing random variable $S>0$ and taking $W(t)= S B(t)- S^2 \sigma^2(t)/2$ in \eqref{Ws}. We retrieve the variance-mixed Brown--Resnick process $\xiW $, which is both max-stable and stationary \citep{eng2012b, Strokorb}.

\subsection{Existence and positivity of $\mathcal{H}_W ^\delta$}
\label{sec2.1}

In order to prove the existence of the generalized Pickands constant $\mathcal{H}_W ^\delta$ we do not need
any further assumptions on the process $W$. In fact, the stationarity of
the process $\xiW$ and the existing theory of max-stable processes
is sufficient to give an immediate answer to a) and partially to b) above.
Indeed, for any compact $E \subset \R$ we define $H_W(E) = \E{ \sup_{t\in E} e^{W(t)}}$
and observe that
\BQN\label{form}
 - \ln \pk{\sup_{t\in E}  \xiW (t) \le x}= H_W(E)e^{-x} , \quad x\inr.
 \EQN
Consequently, by stationarity of $ \xiW $ for any $a\inr$, we have
$H_W(a+E )= H_W(E)$, where $a+E:=\{a+x: x\in E\}.$
Since for any disjoint, non-empty compact sets $E_1, E_2\subset \R$
\begin{align*}
  H_{W}(E_1 \cup E_2)&= \E{\sup_{t\in E_1 \cup E_2} e^{W(t)}}\le  \E{\sup_{t\in E_1} e^{W(t)}}+
  \E{\sup_{t\in E_2} e^{W(t)}}\\
&= H_W(E_1)+ H_W(E_2),
\end{align*}
the set-function $H_W(\cdot)$,
restricted on the sets $\delta \Z \cap [0,T], T>0$, is subadditive and by
Fekete's Lemma
\BQN \label{eq4}
\mathcal{H}_W^\delta=\limit{T} \frac{ H_W(\delta \Z \cap [0,T])}{T}= \inf_{T>0}  \frac{H_W(\delta \Z \cap [0,T])}{T} \in [0,\infty).
\EQN
Therefore, the limit in \eqref{gen_pick} as $T\to\infty$ exists and is finite.
Furthermore, in the case that $\delta >0$,  then \eqref{eq4} immediately implies $\mathcal{H}_W^\delta\le 1/\delta.$

The following lemma is crucial for investigating the structure of $\HWD$
and establishing Dieker--Yakir type representations. It extends Lemma 5.2 in \cite{DM},
where it was considered
for the case that $W(t)=B(t)- \sigma^2(t)/2$ with $B$ a centered Gaussian process with stationary increments and variance function~$\sigma^2$.
\begin{lemma}\label{lemSebastian}
 Suppose that $W$ is such that the process $ \xiW$ in \eqref{eA} is max-stable and stationary, and
 $W(t_0)=0$ for some $t_0\in\R$. If $\Gamma$ is a Borel measurable, positive functional on $D$  that is invariant
under addition of any constant function, then, given that the expectations below exist,
\BQN\label{lemSB}
\E{e^{W(t_0 + t)} \Gamma(W)} = \E{\Gamma(\theta_t W)}, \quad t\inr,
\EQN
where $\theta_t$ is the shift operator, that is, $\theta_t W(s)= W(s-t)$.
 \end{lemma}

An application of equation \eqref{lemSB}  yields a way of rewriting the expectation in
\eqref{gen_pick}; see Corollary 2 in \cite{DiekerY}.

\BEL \label{lemD} \EEH{If  $\mu$ is the Lebesgue measure on $\R$ or the counting measure on  $(k\delta) \Z \cap [0,T]$ with $k\inn, \delta > 0$,  then}
\BQN\label{cor2}
 \frac1T \E{ \sup_{t\in \delta \Z \cap [0,T]} e^{ W (t)}}
&= &     \int_0^1 \mathbb E\left\{\frac{\sup_{s\in \delta\Z \cap [-uT,(1-u)T] } e^{ W (s)}}
   { \int_{-uT}^{(1-u)T} e^{W(s)}\mu(ds) } \right\}\, \mu^T(du) ,
  \EQN
with $\mu^T (du)= \mu ( T du)/T$.
\EEL

\Ec{Using the result of \nelem{lemD}, we establish a Dieker--Yakir representation of $\mathcal{H}_{W}^\delta$
for $\delta>0$} and then show that $\mathcal{H}_{W}^\delta$ is strictly
positive for $\delta\geq 0$.

\BT \label{Propfat} Let $W$ be such that the corresponding max-stable Brown--Resnick
process $\xi_W$ is stationary and has c\`adl\`ag paths. \Ec{If for a given $\delta>0$ we have that
$\pk{S^\delta< \IF}= 1$,} then
\BQN\label{imm}
\mathcal{H}_{W}^\delta = \E{ \frac{ \md} {S^{\delta} }}>0
\EQN
Further, if $\delta \geq 0$ and $\eta= k\delta$ for some $k\inn $, then
\BQN\label{moteta}
 \mathcal{H}_W^\delta \ge
\E{ \frac{ \md }   { S^{\eta} }}>0.
\EQN
\ET

\EEH{The restriction} $\delta >0$ in \eqref{imm} is somehow unsatisfactory. In the sequel we therefore consider
two important special cases where we can strengthen the above results to
\BQN\label{etadelta}
\mathcal{H}_W^{\delta}= \E{ \frac{ \md }   { S^{\eta} }}\in (0,\IF), \quad \delta=0, \eta\ge 0 \text{  or } \delta>0, \eta= k\delta, k\inn,
\EQN
which is motivated by the findings of \cite{DiekerY}  for $W(t)=\sqrt{2}B_\alpha(t)- \abs{t}^\alpha $. Therein \eqref{etadelta}
is shown if $W$ is a fractional Brownian motion and $\delta=0, \eta>0$ or  $\delta= \eta>0$.

\subsection{Gaussian case}
\label{sec2.2}
First, we consider the case where $W(t)=B(t)- \sigma^2(t)/2$, with $B$ a centered, sample continuous Gaussian process that has stationary increments and variance function $\sigma^2$, and
$W(0)=0$ almost surely. In view of \cite{kab2009}, the corresponding $\xiW$ is max-stable and stationary.
In order to apply \netheo{Propfat} we \KD{have to} ensure that
$S^\delta < \IF$ almost surely. To this end, we can require the weak assumption that
\BQN \label{sufMarcus}
\liminf_{\abs{t}\to \IF} \frac{\He{\sigma^2(t)}}{\ln t} > 8,
\EQN
which by Corollary 2.4 in \cite{marcus1972} implies
\BQN\label{IFF}
  \lim_{\abs{t}\to \IF} W(t)= - \IF.
\EQN
Theorem 6.1 in \cite{WangStoev} then yields that $S^\delta < \IF$ almost surely. Consequently,
under \eqref{sufMarcus} and by \netheo{Propfat} we obtain the positivity and Dieker--Yakir representation of  $\mathcal{H}_W^\delta, \delta>0$.

In order to deal with the case $\delta=0$, we need slightly stronger conditions on $\sigma^2$, namely we shall assume that there is an ultimately monotone, non-decreasing function $\ell:[0,\infty)\to [0,\infty)$
and a constant $c\in(0,1]$ such that for all $t$ large
\BQN\label{theta12}
c \ell(t)\le \sigma^2(t) \le \ell(t), \quad  \text{ where   } \limit{t} \frac{ \ell(t)}{\ell(t+k)}=1, \quad \forall k\inn,
\EQN
holds.
Clearly, \eqref{theta12} is satisfied for $\sigma^2$ being a regularly varying function with index $\alpha>0$.
Note in passing that the stationarity of increments implies that $\alpha \le 2$,  \EEH{see also Lemma 2.1 in \cite{marcus1972} for the existence of such Gaussian processes}.
\COM{Additionally, we need to control the speed that $\sigma$ decreases as follows
\BQN\label{theta123}
\sigma^2(t) \le c_* \abs{t}^\beta
\EQN
holds for all $t$ sufficiently small. Condition \eqref{theta123} is satisfied if $\sigma^2$ is regularly varying at 0 with index $\beta>0$. Again, when this is the case, then $\beta \le 2$.}

\BT \label{lb}
Let $W$ be a Gaussian process as above whose variance function $\sigma^2$ satisfies condition \eqref{theta12} with
$c\in (0,1]$ such that $c^2 + 8c -8> 0$. If further
\BQN\label{Sesemn}
 \liminf_{t\to \IF} \frac{\ell(t)}{\ln t}> \frac{8}{c^2 + 8c- 8},
 \EQN
then the generalized Pickands constant $\HWD$ possesses a Dieker--Yakir representation
\BQN
\mathcal{H}_W^{\delta}= \E{ \frac{ M^\delta }   { S^{\eta} }}\in (0,\IF),
\EQN
\Ec{which is valid for $\delta=0$ and $\eta\ge 0$, or $\delta>0$ and $\eta= k\delta, k\inn$.}
\ET
 \BRM
\begin{itemize}
\item[a)]
  Conditions \eqref{theta12} and \eqref{Sesemn}  are much weaker than the assumption that $\sigma^2$ is regularly varying at infinity.  In \cite{debicki2002ruin} the positivity and finiteness of ${\mathcal H}_{W}^0$ is shown under the two conditions C1 and C2 therein, which imply that $\sigma^2$ is a smooth, regularly varying function at infinity and zero.
\item[b)]
  Note that if $c=1$, then \eqref{Sesemn} agrees with \eqref{sufMarcus}.
\item[c)] The validity of \eqref{imm2} can be shown under the assumptions of \netheo{lb}
by borrowing the \KD{arguments of} \cite{DiekerY}.

\end{itemize}
\ERM

\subsection{L\'evy case}
\label{sec2.3}
In \cite{eng2014d}, the so-called L\'evy--Brown--Resnick processes are introduced
as $\xi_W$, where $W$ is composition of two independent L\'evy processes.
More precisely, suppose that $\{B^+(t), t\geq 0\}$ is a L\'evy process such that
its Laplace exponent $\Phi(\theta)=\ln \E{\exp\{\theta B^+(1)\}} $
is finite for $\theta =1$. Define $-W^{-}$ to be the exponentially tilted version of
$$W^{+}(t)=B^+(t)- \Phi(1) t, \quad t\geq 0,$$
that is, the Laplace exponent
of $W^-$ is $\ln \E{\exp\{\theta W^{-}(1)\}}= \Phi(1-\theta)-(1-\theta)\Phi(1)$.
For two independent processes $W^+$ and $W^-$ we
define $W(t)= W^{+}(t), t\ge 0$, and $W(t)=W^{-}(-t)$ if $t< 0$.
With this definition the corresponding process $\xi_W$ is indeed max-stable and
stationary; for details see \cite{eng2014d} and \cite{eng2014c}.

In the case where $B^+$ is a spectrally negative L\'evy process, \cite{eng2014d}
computed the extremal index of the corresponding max-stable process $\xi_W$
explicitly. In view of \eqref{remarkG}, this index coincides with Pickands
constant of the process $B^+$, and it is therefore given as $\KD{\mathcal{H}_W }= \Phi'(1)$.

For more general examples than spectrally negative L\'evy processes,
we show below that
the Pickands constant $\mathcal{H}_W$ in the L\'evy case
possesses a Dieker--Yakir type representation.
In fact, by \cite{eng2014d} it follows that the conditions of \netheo{Propfat} are satisfied
and thus $\HWD$ exists and is strictly positive.
In what follows we suppose that $B^+$ is not a Poisson process with
lattice support of jump distribution.

\BT \label{levy}
Let $B^+(t),t\in [0,\IF)$ and $W(t),t\inr$ be as above. \\
(1) If $\E {e^{(2+\varepsilon)|W(1)|}}<\infty$
and
$\E {e^{(2+\varepsilon)|W(-1)|}}<\infty$
for some $\ve>0$, then
\BQN
\mathcal{H}_W^0= \E{ \frac{ M^0 }   { S^0 }}\in (0,\IF). 
\EQN
(2) If $\E {e^{(1+\varepsilon)|W(1)|}}<\infty$
and
$\E {e^{(1+\varepsilon)|W(-1)|}}<\infty$
for some $\ve>0$, then
\BQN
\mathcal{H}_W^\delta= \E{ \frac{ M^\delta }{ S^\eta }}\in (0,\IF),
\quad \delta=0, \eta> 0 \text{  or } \delta>0, \eta= k\delta, k\inn. 
\EQN
\ET

\KD{
\BRM
Theorem \ref{levy} holds if both the left and the right tail
probability of $W(1)$ is sufficiently light; for example if
$\Phi(\theta)<\infty$ for $\theta\in (-2-\varepsilon,3+\varepsilon)$ for scenario (1)
and  $\theta\in (-1-\varepsilon,2+\varepsilon)$ for scenario (2).
We conjecture that the claim of Theorem \ref{levy}
is true under weaker assumptions on $W$.
\ERM
}

\section{A connection to mixed moving maxima processes}\label{sec_m3}
As in the previous section, let $W$ with $W(0)=0$ a.s.\ be a c\`adl\`ag process
such that the corresponding $ \xiW $ is max-stable and stationary.
The process $ \xiW $ is said to admit
a mixed moving maxima representation (for short M3) if
\begin{align}\label{M3}
   \xiW (t) \stackrel{d}{=} \max_{i\ge 1} (F_i(t- P_i)+ Q_i), \quad t\inr,
\end{align}
where the $F_i$'s are independent copies of a measurable c\`adl\`ag
process $F_W(t), t\inr$, with \BQN\label{normaz} \sup_{t\in \R} F_W(t)
= F_W(0) = 0 \EQN almost surely, and \BQN C_W= \left(\E{ \int_{\R}
  \exp( F_W(t)) \, dt}\right)^{-1} \in (0,\IF).  \EQN Here,
$\sum_{i=1}^\IF \ve_{(P_i,Q_i)}$ is a Poisson point process in $\R^2$
with intensity $C_W dt \, \exp(-y)dy$.  Note that the normalization of
the supremum of $F_W$ to $0$ in \eqref{normaz} is crucial since
otherwise the constant $C_W>0$ would not be well-defined. Furthermore,
$C_W$ ensures that the margins of $ \xiW $ are standard Gumbel
distributions and it appears thus naturally in the theory of
max-stable processes. It plays a crucial role in the simulation of
such processes but its numerical evaluation is time intensive and the
exact value is, apart from special cases, unknown  \citep{KabExt}.

Throughout this section we assume that $\xi_W$ possesses a M3 representation
which amounts to assuming one of the equivalent conditions below; for details see
\cite{WangStoev} and Theorem 2 in \cite{dom2016}.
\begin{cond}\label{cond_mmm}
  We assume that one of the following equivalent conditions holds:
  \begin{enumerate}
  \item
    The max-stable process $ \{\xiW, t\inr\}$ possesses a M3 representation.
  \item
    The max-stable process $ \{\xiW, t\inr\}$ has no conservative component in its spectral representation.
  \item
    The process $\{W(t), t\inr\}$ satisfies
    \begin{align*}
      \lim_{|t|\to\infty} W(t) = -\infty, \quad \text{a.s.}
    \end{align*}
  \item
    The process $\{W(t), t\inr\}$ fulfills
    \begin{align*}
      \int_{-\infty}^\infty  e^{W(t)} dt < \infty, \quad \text{a.s.}
    \end{align*}
  \end{enumerate}
\end{cond}

Since we are interested also in the case $\delta >0$, we \BH{show}
in the Appendix how to derive an M3 representation for the discretized process
$ \xiW ^\delta= \{\xiW (t), t\in \delta \Z\}$, with shape functions $F_W^\delta$ and constant
\begin{align*}
  C^\delta_W = \left(\E{ \int_{t\in \delta\Z} \exp( F_W^\delta(t)) \nu_\delta(dt)}\right)^{-1} \in (0,\IF), \quad \delta > 0.
\end{align*}
Here, $\nu_\delta / \delta$ for $\delta > 0$ is the counting measure on $\delta \Z$,
and $\nu = \nu_0$ is the Lebesgue measure. In the sequel the superscript is omitted if it is 0, for instance we write
$C_W$ and $\mathcal{H}_W$ instead of $C_W^0$ and $\mathcal{H}_W^0$, respectively.
The M3 representation of $\xiW$ allows us to show a new formula for $\mathcal{H}_W^\delta$ and relate it to $C_W^\delta$.
Moreover, we prove that  $C_W^\delta$ is exactly what we refer to as the Dieker--Yakir representation of Pickands constant.

\begin{theorem} \label{theothree}
\EEH{If $ \xiW^\delta $ possesses an M3 representation, then for any $\delta\ge 0$}
\BQN\label{posH}
0 <   C_W^\delta = \E{\frac{ \md}{S^\delta }} \le \He{\mathcal H_W^\delta}.
\EQN
 \end{theorem}

\BRM For any $\delta > 0$, in view of
\eqref{imm} in Theorem \ref{Propfat} and \eqref{posH}, we have the equality
$$ \HWD= C_W^\delta>0.$$
\ERM

 Except \KD{for} few special cases, the exact value of Pickands constant
is unknown. There are several attempts
to assess its value by Monte Carlo simulation, most notably via the
recent Dieker--Yakir representation in \cite{DiekerY}.
The above Corollary states that the simulation problem of Pickands constant
$\mathcal H_W^\delta$ is equivalent to the problem of simulating the
constants $C_W^\delta$ in spatial extreme value theory, provided that $ \xiW $
admits an M3 representation and the Dieker--Yakir representation for $\mathcal{H}_W^\delta$ holds. This is a fruitful observation since there is active research
on the simulation of max-stable processes \citep{DM, dom2016a} and even
of the constant $C_W^\delta$ \citep{KabExt}.  
We conclude this section with several examples.

\begin{example}
If $W(t) = \sqrt 2 Z t - t^2$, $t\in\R$, where $Z$ is
an $N(0,1)$ random variable it is known \citep{WangStoev} that $\xiW$
has an M3 representation with
deterministic shape functions
 $   F_W(t) = - t^2, t\in \R.$
Thus
  \begin{align*}
C_W=\Bigl( \int_\R e^{- t^2}\,  dt\Bigr)^{-1}  = \frac {1}{\sqrt \pi},
  \end{align*}
and consequently, by \netheo{lb} and \ref{theothree} we recover the well-known fact $ \mathcal H_{W} =1/\sqrt{\pi}$.

If $W(t) =\sqrt 2 B(t)- \abs{t},t\inr$,  where $B$ is
a standard Brownian motion, then it follows by \citet{eng2011} that
$\xiW$ has an M3 representation whose shape functions $F$ are
given by a three-dimensional Bessel process and that $ C_W = 1.$
Thus, again by \netheo{lb} and \ref{theothree} we recover
 $\mathcal{H}_W=C_W= 1$ \citep{Pit20}. 
\end{example}

\begin{example} \label{rem34} Suppose that $W$ is a sample continuous Gaussian process with stationary
  increments that fulfills the assumptions of Theorem \ref{lb}.
  Since in this case \eqref{sufMarcus} holds and thus Condition \ref{cond_mmm}
  is satisfied,
  $ \xiW $ admits an M3 representation and
  in view of  \netheo{theothree}
  $\mathcal H_W^\delta$ is positive  for any $\delta>0$ and 
  $$C_W^\delta= \frac{1}{\delta} \pk{ \sup_{t\in \delta \Z} W(t)=0}= \HWD.$$
Furthermore, we have
  $$\lim_{\delta\downarrow 0} C_W^\delta=C_W= \mathcal{H}_W =\E{\frac{ M^0}{S^0 }}.$$
\end{example}

\begin{example}\label{remB}
If $W$ is as in Section \ref{sec2.3}, \cite{eng2014c} show that the L\'evy-Brown--Resnick process $\xiW$ admits an M3 representation  where the constant   $C_W$ is explicitly given by
$$C_W = \frac{\underline{k}(0,1)}{\underline{k}'(0,0)} > 0,$$
where $\underline k$ is the bivariate Laplace exponent of the descending
ladder process corresponding to $W$.
In particular, this implies that for $\delta=0$ by \netheo{theothree} $\mathcal{H}_W\geq C_W $ and thus
\BQN
\mathcal{H}_W \ge \frac{\underline{k}(0,1)}{\underline{k}'(0,0)} > 0.
\EQN
In order to have equality in the equation above, it is sufficient that the process $W$ satisfies
the conditions of \netheo{levy}, since then
$$\mathcal{H}_W   = \E{\frac{ M^0}{S^0 }} = C_W.$$
\end{example}

\section{Proofs}\label{proofs}

\prooflem{lemSebastian} It is well-known that the stationarity of $ \xiW $ is equivalent to the fact
  that for arbitrary $h \in \R$ the two Poisson point processes
  $\{ U_i + W_i: i\in\N\}$
  and $\{ U_i + \theta_hW_i: i\in\N\}$ on $D$  have the same
  intensity; see \cite{kab2009}. The latter holds if and only if for any Borel subset $A$ 
   \begin{align*}
    \int_{\R} e^{-u} \pk{ u + W \in A }du = \int_{\R} e^{-u} \pk{ u + \theta_h W \in A }du.
  \end{align*}
   Let $B \subset D$ be a shift-invariant Borel set in the sense that  $B + x = B$ for any $x\in \R$,
   and recall that $W(t_0)=0$ almost surely.
   Consequently, for any $h\inr$ we have
   \BQNY
   \E{e^{W(t_0 +h)} \mathbf 1\{W\in B\}} &= &\E{\int_{\R} e^{-u} \mathbf 1\{u + W(t_0+h) > 0 \}\mathbf 1\{W\in B\} du }\\
   & =& \int_{\R} e^{-u} \pk{ u + W(t_0+h) > 0, u + W\in B} du \\
   & =& \int_{\R} e^{-u} \pk{ u + W(t_0) > 0, u + \theta_h W\in B}  du\\
   & = &\int_{\R} e^{-u} \mathbf 1\{u  > 0 \}\pk{ u + \theta_h W\in B} du\\
   & = &\pk{ \theta_h W\in B}.
   \EQNY
   Furthermore, the above readily extends to Borel measurable, positive functionals $\Gamma$ on $D$ that are invariant under addition of a constant function and, thus, the assertion follows. \QED

\prooflem{lemD} Define the translation invariant functional
  \begin{align*}
    \Gamma(f) = \frac{\sup_{s\in \delta \Z \cap [0,T]} e^{ f (s)}}
   {\int_{0}^{T} e^{f(s)}\mu(ds )}.
  \end{align*}
Clearly, we have that for any $t\in (k\delta) \Z$
\BQNY
 \Gamma(\theta_t f) &=&
\frac{\sup_{s\in \delta \Z \cap [0,T]} e^{ f (s-t)}}
   {\int_{0}^{T} e^{f(s-t)}\mu(ds )} =
\frac{\sup_{s\in \delta \Z \cap [-t,T-t]} e^{ f (s)}}
   {\int_{0}^{T} e^{f(s-t)}\mu(ds )} = \frac{\sup_{s\in \delta \Z \cap [-t,T-t]} e^{ f (s)}}
   {\int_{-t}^{T-t} e^{f(s)}\mu(ds )},
   \EQNY
   where the last equality follows by the translation invariance of $\mu$.
Hence, as in the proof of Corollary 2 in \cite{DiekerY} a direct  application of \nelem{lemSebastian} yields
  \BQNY
 \frac1T \E{ \sup_{t\in \delta \Z \cap [0,T]} e^{ W (t)}}
&=& 
    \frac1T \int_0^T \mathbb E\left\{e^{W(t)} \frac{\sup_{s\in \delta \Z \cap [0,T] } e^{ W (s)}}
   {\int_{0}^{T} e^{W(s)}\mu(ds ) } \right\} \mu(dt) \\
&=&     \frac1T \int_0^T  \E{ e^{W(t)}\Gamma(W)}\mu(dt)\\
&=&     \frac1T \int_0^T  \E{ \Gamma(\theta_t W)}\mu(dt)\\
   &=&\frac1T \int_0^T \mathbb E\left\{ \frac{\sup_{s\in \delta\Z \cap [-t,T-t] } e^{ W (s)}}
   {\int_{-t}^{T-t} e^{W(s)}\mu(ds) } \right\} \mu(dt).
   \EQNY
   Consequently, \eqref{cor2} follows by changing the variable $t = u T$. \QED

\prooftheo{Propfat}
Let first $\eta=\delta>0$, then if $\BH{\lambda_\delta}$ denotes the counting measure on $\delta \Z$,
then applying \eqref{cor2} with  $\mu=\lambda_\delta$ we obtain
$$\frac1T \E{ \sup_{t\in \delta \Z \cap [0,T]} e^{ W (t)}} =
    \int_0^1 \mathbb E\left\{\frac{\sup_{s\in \delta\Z \cap [-uT,(1-u)T] } e^{ W (s)}}
   { \delta \int_{-uT}^{(1-u)T} e^{W(s)}\mu(ds) }  \right\}\,  \delta \mu^T(du). $$
By the assumption that $S^\delta= \delta \int_{\R} e^{W(s)}\mu(ds) < \IF$ it follows that
$\sup_{s \in \delta \Z} e^{W(s)} < \IF$ and $\lim_{\abs{n}\to \infty, n\in \Z} W(n \delta ) = -\infty$ almost surely.
Hence the almost sure convergence
$$    g_{T,\delta}(u)= \frac{\sup_{s\in \delta\Z \cap [-uT,(1-u)T] } e^{ W (s)}}
   {\delta \int_{-uT}^{(1-u)T} e^{W(s)}\mu(ds) }   \to \frac{\sup_{s\in \delta\Z  } e^{ W (s)}}
   {\delta\int_{\R} e^{W(s)}\mu(ds) }  = \mathcal{Q}_\delta \le \frac{1}{\delta}, \quad T \to \IF $$
holds    for any $u\in (0,1)$. Clearly, the above convergence remains true if we
replace $u$ by a sequence $u_T, T>0$ such that $\limit{T} u_T=  u \in (0,1)$.
Since for any $u\in(0,1),T>0$ we have $g_{T,\delta}(u)\le 1 / \delta$ we obtain for any $u\in(0,1)$
by dominated convergence
$$ \limit{T} \E{g_{T,\delta}(u_T)}= \E{\mathcal{Q}_\delta}.$$
Since $\delta \mu^T$ converges weakly to the Lebesgue measure as $T\to \IF$,
Theorem 5.5 in \cite{bil1968} implies that
$$\HWD =  \limit{T} \int_0^1 \E{g_{T,\delta}(u)} \delta \mu^T(du) = \int_0^1 \E{\mathcal{Q}_\delta} \, du = \E{\mathcal{Q}_\delta} $$
establishing the first claim in \eqref{imm}. \\
Next, if $\mu=\BH{\lambda_\eta}$ with $\eta=k\delta, k=0,1 \ldot$ or $\eta>0, \delta=0$, \BH{by  \eqref{cor2} and }
Theorem 1.1 in \cite{Fatou14} for any $u\in (0,1), T>0$
\BQNY
\mathcal{H}_W^\delta&= & \limit{T} \int_0^1\E{g_{T,\eta}(u)  }  \nu^T_\eta(du)  \\
&\ge &  \int_0^1 \liminf_{T\to \IF, v\to u} \E{g_{T,\eta}(v)} \, du  \\
&\ge &  \int_0^1 \E{\liminf_{T\to \IF, v\to u} g_{T,\eta}(v)}\, du \\
&= &\E{\mathcal{Q}_\eta}>0,
\EQNY
hence \eqref{moteta} follows. \QED

\prooftheo{lb} Our assumptions on $\sigma^2$ imply that \eqref{sufMarcus} holds, and thus $\delta \sum_{t\in \delta \Z} e^{W(t)}< \IF$ almost surely for any $\delta\ge 0$.
Recall that we interpret $\delta \sum_{t\in \delta Z} e^{W(t)}$ as  $\int_{\R} e^{W(t)}\, dt$ when $\delta=0$.
  Consequently, for any $\delta,\eta\ge 0$, we have the almost sure convergence
\begin{align}\label{asRatio}
  R_{u,T}^{\delta,\eta}&=\frac{M^\delta{[-uT, (1-u)T]}} {S^\eta{[-uT, (1-u)T]} } \\
  \notag& := \frac{\sup_{s\in \delta \Z \cap [-uT, (1-u)T]}e^ {W(s)}}{\eta \sum_{s\in \eta \Z  \cap [-uT, (1-u)T]} e^{W(s)}}\to
\frac{\sup_{s\in \delta \Z} e^{W(s)}} { \eta \sum_{s\in \eta \Z  } e^{W(s)}} \in (0, \IF)
\end{align}
for all $u\in (0,1), T\to \IF$.
Together with \eqref{cor2}, the claim of the theorem therefore follows if we can show the uniform integrability
$$ \limit{A} \sup_{T>0} \sup_{u\in (0,1)} \E{R_{u,T}^{\delta,\eta}; R_{u,T}^{\delta,\eta}> A  } = 0.$$
In order to give a self-contained proof (which follows along the same ideas as in \cite{DiekerY})
we introduce the same notation as therein. Namely, we
let $a_j=j$ and we define
$J_j=[a_j,a_{j+1})$,
$J_{-j}=(-a_{j+1},-a_{j}]$,
$S_j^\eta=\eta \sum_{k:\eta k\in J_j} e^{W(\eta k)}$,
$M^\eta_j=
\sup_{k: \eta k\in J_j} e^{W (\eta k)}
$
and
$S^\eta,  M^\eta$ as in (\ref{md}).
Note that in the aforementioned paper our $W$ corresponds to $Z$.

Fix some $\lambda>0$ and define  $W^\lambda(t)= W(\lambda \floor[\big]{t/\lambda}), t>0$, and $W^\lambda(t)= W(\lambda \ceil[\big]{t/\lambda})$
otherwise. We have
$$ \frac{M_j^\lambda}{S_j^\lambda}\le \frac{1}{\lambda}, \quad
\frac{M_j^\delta}{ M_j^\lambda}  = e^{\sup_{s\in J_j} W^\delta(s)-\sup_{s\in J_j} W^\lambda(s)} \le
 e^{\sup_{s\in J_j} \abs{W^\delta(s)-W^\lambda(s)}}$$
 and
$$   \frac{ S_j^\lambda}{S_j^\eta}
= \frac{\int_{a_j}^{a_{j+1}} e^{W^\lambda(t)- W^\eta(t)}e^{W^\eta (t)}\, dt} {\int_{a_j}^{a_{j+1}} e^{W^\eta (t)}\, dt}  \le
e^{\sup_{s\in J_j}\abs{ W^\eta (s)-W^\lambda(s)}}.$$

On the event $\{M^\delta =M_j^\delta\}$ for some $j\in \Z$ we have
(assume that $uT, (1-u)T\in \Z$)
\BQN \label{00}
R_{u,T}^{\delta,\eta} &\le & \frac{M_j^\delta}{ S_j^\eta} = \frac{M_j^\lambda}{S_j^\lambda}\frac{M_j^\delta}{ M_j^\lambda}  \frac{ S_j^\lambda}{S_j^\eta} \notag \\
&\le & \frac{1}{\lambda} e^{ \sup_{s\in J_j}\abs{B^\delta(s)- B^\lambda(s)} +
\sup_{s\in J_j}\abs{B^\eta(s)- B^\lambda(s)}+ \kappa_\lambda(j)}=:R_{u,T}^{\delta,\eta}(j), 
\EQN
where $B(t)= W(t)+\sigma^2(t)/2$ is a centered Gaussian process and
$$ \kappa_\lambda(j):= \sup_{s\in J_j} \abs{Var(W(s))- Var(W^\lambda(s))}. $$
Since $M^\delta{[-uT, (1-u)T]} \ge 1$, we have
\BQNY
\lefteqn{\E{R_{u,T}^{\delta,\eta}; R_{u,T}^{\delta,\eta}> A  }}\\
 &=&
\sum_{j \in \Z} \E{R_{u,T}^{\delta,\eta}; R_{u,T}^{\delta,\eta}> A, M_j^\delta=M^\delta  } \\
&\le & \E{R_{u,T}^{\delta,\eta}; R_{u,T}^{\delta,\eta}> A, M_0^\delta=M^\delta  }+ 2
\sum_{j\ge 1} \E{R_{u,T}^{\delta,\eta}(j); R_{u,T}^{\delta,\eta}(j)> A, M_j^\delta\ge 1  }  \\
&\le & \E{M_0^\delta/\EEH{S_0^\eta};M_0^\delta/ \EEH{S_0^\eta} >A } \\
&& \hspace{2cm} + 2\sum_{j\ge 1} \E{R_{u,T}^{\delta,\eta}(j); R_{u,T}^{\delta,\eta}(j)> A, \sup_{s\in J_j} B^\delta(s)\ge \inf_{s\in J_j} \sigma^2(\delta\floor{s/\delta})/2  }  \\
&=:& \E{M_0^\delta/\EEH{S_0^\eta};M_0^\delta/\EEH{S_0^\eta} >A }+ 2 \sum_{j\ge 1} \pi_{j}(A).
\EQNY
In the following $C>0$ may change from line to line.
We note that 
$$\E{M_0^\delta/\EEH{S_0^\eta};M_0^\delta/\EEH{S_0^\eta} >A } \to 0$$
as $A\to \IF$ since $ \E{M_0^\delta/\EEH{S_0^\eta}} < \infty$.
For all $t,s \in J_j$  and by \eqref{theta12} for all  $a_j$ large enough, by the monotonicity of $\ell$
$$  \inf_{s\in J_j} \sigma^2(\delta \floor{s/\delta})/2 \ge c \ell(\delta \floor{a_j/\delta}),  \quad
\sup_{s\in J_j} Var(B^\delta(s)) \le  \ell(\delta \floor{(a_{j}+1)/\delta}) .$$
Since for all $j\in \Z$
\BQNY
\E {\sup_{s\in J_j} B^\delta (s)} &=&
\E {\sup_{s\in J_j} [B^\delta (s)- B^\delta(j)+ B^\delta(j)]}\\
&=& \E {\sup_{s\in [0,1]} (B^\delta (s+j)-B^\delta (j))}\\
&\le&  C
\EQNY
where the last inequality is consequence of
$$ \sup_{s\in [0,1]}Var(B^\delta(s+j) - B^\delta(j)) = \sup_{s\in [0,1]}
\sigma^2(\delta [\floor{(s+j)/\delta}- \floor{j/\delta}])< C,$$
then by Borell-TIS inequality (see, e.g., \cite{GennaBorell})
\BQNY
\pk{ \sup_{s\in J_j} B^\delta (s)\ge \inf_{s\in J_j} \sigma^2(\delta \floor{s/\delta})/2 } & \le &
\pk{ \sup_{s\in J_j} B^\delta (s)\ge c\ell( \delta \floor{a_j/\delta})/2 }\\
&\le & 
C\exp\left((1- \ve_1) \frac{c^2\ell^2(\delta \floor{a_j/\delta})}{8 \ell(\delta \floor{(a_{j}+1)/\delta} )} \right) \\
&\le & C \exp\left( - (1- \ve_2) \frac{c^2 \ell(a_j+1)}{8}  \right)
\EQNY
for some $\ve_1,\ve_2$  positive arbitrary small \EEH{and all} $j\ge 1$. Further, the fact that
$$ \sup_{s\in J_j}Var(B^\delta(s) - B^\lambda(s)) = \sup_{s\in J_j}
\sigma^2(\delta \floor{s/\delta}- \lambda \floor{s/\lambda})< C$$
for \EEH{all} $j$, that is, the variance is bounded implies (use Borell-TIS inequality, see e.g., \cite{AdlerTaylor})
$$ \E{ e^{ p \sup_{s\in J_j} \abs{B^\delta(s)- B^\lambda(s)}  }}\le C$$
for any $p> 1$ \EEH{and all} $j$.  Consequently, by the H\"older inequality for $q=1+ 1/(p-1)$ and $\ve>0$ sufficiently small
\BQNY
\pi_j(A) &\le & \frac1\lambda e^{\kappa_\lambda(j)}\Biggl(\E{ e^{
p \sup_{s\in J_j} \abs{B^\delta(s)- B^\lambda(s)} +
p \sup_{s\in J_j} \abs{B^\eta(s)- B^\lambda(s)}  }} \Biggr)^{1/p} \\
&& \times \Biggl( \pk{R_{u,T}^{\delta,\eta}(j)> A }\Biggr)^{1/(pq)}\Biggl( \pk{ \sup_{s\in J_j} B^\delta(s)\ge \inf_{s\in J_j} \sigma^2(\delta \floor{s/\delta})/2 } \Biggr)^{1/q^2}\\
&\le & C \Biggl( \pk{R_{u,T}^{\delta,\eta}(j)> A }\Biggr)^{1/(pq)} \exp\left(\kappa_\lambda(j) - (1- \ve_2) \frac{c^2 \ell(a_j+1)}{8q^2} \right).
\EQNY
Further, by our assumptions on $\ell$ and $c$, for all $j$ large and $\ve_3>0$ sufficiently small
\BQNY
\lefteqn{\kappa_\lambda(j)= \sup_{s\in J_j} \abs{ \sigma^2(s)- \sigma^2(\lambda \floor{s/\lambda})} }\\
&\le &
\max\Biggl( \ell(a_j+1)   -  c \ell(\lambda \floor{a_j/\lambda})   ,
 \ell(\lambda \floor{(a_j+1)/\lambda}) -  c \ell(a_j)  \Biggr)  \\
&\le& \EEH{(1- c+\ve_3)} \ell(a_j+1).
\EQNY
Choose $q>1$ sufficiently close to 1. \KD{Then, by the assumption $c^2 + 8c -8> 0$ and in view of \eqref{Sesemn},}
we can find a constant $B> 1$ and take $\ve_i>0, i\le 4 $, sufficiently small  such that
\begin{align*}
\sum_{j\ge 1} &\exp\left(\kappa_\lambda(j)- (1-\ve_2)\frac{c^2}{8q^2} \ell(a_j+1)  \right)\\
&\le 
\sum_{j\ge 1}  \exp\left(-  \Bigl((1- \ve_2) \frac{c^2}{8 q^2} - (1-c+\ve_3)\Bigr) \ell(a_j+1)  \right)\\
&\le
\sum_{j\ge 1}  \exp\left(- \EEH{ (c^2 + 8c -8- \ve_4)} \ell(a_j)   \right)\\
&\le
\sum_{j\ge 1}e^{- B \ln a_j}= \sum_{j\ge 1}\frac{1}{a_j^{B}} < \IF.
\end{align*}
Therefore
\BQNY
\sum_{j\geq 1} \pi_j(A) \leq  \sum_{j\geq 1} C \Biggl( \pk{R_{u,T}^{\delta,\eta}(j)> A }\Biggr)^{1/(pq)} \frac{1}{a_j^{B}} \to 0, \quad A\to\infty,
\EQNY
which concludes the proof. \QED

\prooftheo{levy}
The idea of the proof is similar to the proof
of Theorem \ref{lb}, with slight modifications
which we analyze below. We use the same notation as in the proof of the aforementioned theorem and focus on the case that $uT, (1-u)T\in \Z$.
\\
{\underline{\it Case $\eta=0$}}.
Since $\delta=0$ in this case, we set $\lambda=1$ and observe that,
on the event $\{M^\delta =M_j^\delta\}$,
\begin{eqnarray*}
R_{u,T}^{0,0}
&\le &
e^{ 2\sup_{s\in J_j}\abs{W(s)- W^1(s)}}
=:\widehat{R}_{u,T}^{0,0}(j)
\end{eqnarray*}
and
\BQNY
\E{R_{u,T}^{0,0}; R_{u,T}^{0,0}> A  }
&\le &
\E{R_{u,T}^{0,0}; R_{u,T}^{0,0}> A, M_0^0=M^0  }\\
&& +
\sum_{j\in \Z\setminus\{0\}} \E{\widehat{R}_{u,T}^{0,0}(j);
\widehat{R}_{u,T}^{0,0}(j)> A, M_j^0\ge 1  }\\
&=:&
\E{R_{u,T}^{0,0}; R_{u,T}^{0,0}> A, M_0^0=M^0  }
+
\sum_{j\in \Z\setminus\{0\}}
 \widehat{\pi}_{j}(A).
\EQNY

As in the proof of Theorem \ref{lb},
$\limit{A}\E{M_0^0/S_0^0;M_0^0/S_0^0 >A } = 0$
since $ \E{M_0^0/S_0^0 } < \infty$.
Thus we focus on an upper bound for
$ \widehat{\pi}_{j}(A)$.
By the same argument as given in the proof of Theorem \ref{lb},
for any $p>1$ and $q=1+ 1/(p-1)$,
\BQNY
\widehat{\pi}_j(A) &\le &  \Biggl(\E{ e^{
2p \sup_{s\in J_j} \abs{W(s)- W^1(s)} }} \Biggr)^{1/p}
\Biggl( \pk{\widehat{R}_{u,T}^{0,0}(j)> A }\Biggr)^{1/(pq)}
          \Biggl( \pk{ \sup_{s\in J_j} W(s)\ge 1}  \Biggr)^{1/q^2}.
\EQNY
Suppose that $j \ge 1$.
By (2.1)  in \cite{Wil87} (see also Lemma 9.1 in \cite{DeM15}), for each $u>u_0>0$
\[
\pk{ \sup_{s\in [0,1) } W(s)>u }\le \frac{\pk{W(1)>u-u_0}}{\pk{\inf_{s\in [0,1) } W(s)>-u_0}}
\]
and
\[
\pk{ \inf_{s\in [0,1) } W(s)<-u }\le \frac{\pk{-W(1)>u-u_0}}{\pk{\inf_{s\in [0,1) } -W(s)>-u_0}},
\]
which implies that
\begin{align}
\label{lev.1} \E{ e^{2p \sup_{s\in J_j} \abs{W(s)- W^1(s)} }}
&=\E{ e^{2p \sup_{s\in [0,1) } \abs{W(s)} }}\\
\notag &\le C_1 \E{ e^{ 2p \abs{W(1)} }}<\infty
\end{align}
for sufficiently small $p>1$ and some $C>0$.

Next, in order to derive a tight upper bound for
$\pk{\sup_{t\in J_j}W(t)>1}$, as $j\to\infty$, let us recall that $W(t)=B^+(t)-\Phi(1)t$, for $t\ge0$,
and observe that
$\E {W(1)} = \E{B^+(1)-\Phi(1)}<0$.

Let $\varepsilon=\frac{1}{2}(\Phi(1)-\E {B^+(1)})>0$
and introduce the following L\'evy process
$L(t):=W(t)+\varepsilon t$.
It is straightforward to check that
$\E {L(1)}<0$ and for
$\Phi_L(\theta):=\ln \E {e^{\theta L(1) }}   $
we have
$$\Phi_L(0)=0, \quad \Phi'_L(0)=\E {L(1)}<0$$
and $\Phi_L(1)=\varepsilon>0$. Hence, there exists $1>\gamma>0$ such that
$\Phi_L(\gamma)=0$ and
$\Phi'_L(\gamma)<\infty$.
Now, following, e.g., Theorem 2.6 from \cite{AsA10}
\begin{eqnarray}
\hspace{2em} \pk{\sup_{t\in J_j}W(t)>1 }
\le
\pk{  \sup_{t\in[0,\infty)}(B^+(t)-t(\Phi(1)-\varepsilon))>\varepsilon j}
\le {\rm C} e^{-\gamma \varepsilon j} \label{lev.2}
\end{eqnarray}
for some ${\rm C}\in (0,\infty)$ and all $j\ge 1$. Therefore, combining (\ref{lev.1}) with (\ref{lev.2}), we get
\BQNY
\sum_{j\geq 1} \widehat{\pi}_j(A)\to 0, \ \ \rm{as} \ \ A\to\infty.
\EQNY
The proof that $\limit{A}\sum_{j\leq -1} \widehat{\pi}_j(A)=0$ follows by the same argument, with the use of the fact that
$W(t)=W^{-}(-t)$ if $t< 0$, with
$\ln \E{e^{\theta W^{-}(1)}}= \Phi(1-\theta)-(1-\theta)\Phi(1)$.

{\underline{\it Case $\eta>0$}}.
We set
$a_j:=\eta j$
and
observe that, on the event $\{M^\delta =M_j^\delta\}$ (assume that $uT, (1-u)T\in \Z$),
\BQN
R_{u,T}^{\delta,\eta}
&\le &
\frac{M_j^\delta}{ S_j^\eta}
=
\frac{M_j^\delta}{M_j^\eta}  \frac{ M_j^\eta}{S_j^\eta} \notag
\le
\frac{1}{\eta}
e^{\sup_{s\in J_j} \abs{W^\delta(s)-W^\eta(s)}}
\le
\frac{1}{\eta}
e^{\sup_{s\in J_j} \abs{W(s)-W^\eta(s)}}.
\EQN
The rest of the proof goes line by line the same as the proof of case $\eta=0$,
with the use of the fact that if $\eta>1$, then
$\E{ e^{  p|W(\eta)|   }   }\le \left(\E{ e^{  p|W(1)|   }   }\right)^{\lceil \eta \rceil}$.
This completes the proof.

\QED

\prooftheo{theothree}
 For an M3 process as above, the finite dimensional distributions of $ \xiW ^\delta$
for $t_i\in \delta \Z,x_i\inr, 1\leq i\le n$, $n\inn$ can be written as
\begin{align}\label{crucFw}
- \ln &\pk{\xi^\delta_ W(t_1) \le x_1 \ldot \xi^\delta_ W(t_n)\le x_n}\\
   \notag&\hspace{2.5cm}=   C_W^\delta \E{ \int_{r\inr}
  \max_{j=1,\dots, n} \exp\left(F_W^\delta(t_j-r)- x_j\right)
\, \nu_\delta(dr)}.
\end{align}
Since $ \xiW $ has c\`adl\`ag paths by assumption, we have for any compact set $E\subset \R$
\BQNY
- \ln \pk{ \sup_{t\in  \delta \Z \cap E}  \xiW ^\delta(t) \le \He{0} }&=&   C_W^\delta \E{ \int_{r\in \delta \Z} \sup_{t\in
 \delta \Z \cap E} \exp\left( F_W^\delta(t-r) \right)\, \nu_\delta(dr)}, 
\EQNY
which, in view of equation \eqref{remarkG},  implies
\BQN\label{Aa}
\He{\mathcal H_W^\delta} =  \limit{T} \frac{C_W^\delta}{T}  \E{\int_{r\in \delta \Z} \sup_{t\in \delta \Z \cap [0,T]} \exp\left( F_W^\delta(t-r) \right)\, \nu_\delta(dr) }.
\EQN
Set $T_\delta=T$ if $\delta=0$ and $T_\delta= \delta \lfloor T/\delta \rfloor$ otherwise. For any fixed $T>0$
\begin{align*}
  &\E{ \int_{r\inr} \sup_{t\in \delta \Z \cap [0,T]} \exp\left(F_W^\delta(t-r) \right)\, \nu_\delta(dr)} \\
  &\hspace{3cm}\geq \E{ \int_{-T_\delta}^0 \sup_{t\in \delta \Z \cap [0,T]}\exp\left(F_W^\delta(t-r) \right)\, \nu_\delta(dr)}
\, dr = T_\delta,
\end{align*}
 since by the assumption $\sup_{t\in \delta \Z}  F_W^\delta(t) = F_W^\delta(0) = 0$ almost surely, for any  $r\in[- T_\delta,0]$ we have
 $$\sup_{t\in \delta \Z \cap [0,T]} \exp\left(F_W^\delta(t-r) \right) = \exp\left(F_W^\delta(0) \right)=1.$$
Consequently,
$$ \mathcal{H}_W^\delta \ge C_W^\delta\limit{T} \frac{T_\delta}{T}= C_W^\delta.$$
We show next \eqref{posH} for  $\delta = 0$.  Theorem 4.1 in \cite{engelke14} implies that
  the process $W$ with $W(0)=0$ almost surely, can be obtained by the M3 representation in terms of the shape
  function as
  \begin{align*}
    \pk{W\in L} = C_W \int_{D} \int_\R \mathbf{1} \left\{ f(\cdot + s) - f(s) \in L\right\} e^{ f(s)} ds \mathbb P_{F_W}(df),
  \end{align*}
  where $L$ is an arbitrary Borel subset of $D$. Consequently,  for any $\mathbb P_W$-integrable functional $\Gamma: D \to \R$ we have
  \begin{align*}
\E{\Gamma(W)} = C_W \int_{D} \int_\R \Gamma(f(\cdot + s) - f(s)) e^{f(s)} ds \mathbb P_{F_W}(df).
  \end{align*}
  Let now $\Gamma$ be given by the mapping (on a suitable subspace of $D$  with full $\mathbb P_W$ measure)
  \begin{align*}
    f \mapsto \frac{\sup_{t\in\R} e^{f(t)}}{\int_\R e^{f(t)} dt},
  \end{align*}
  and observe
  \begin{align*}
\E{ \frac{\sup_{t\in\R} e^{W(t)}}{\int_\R e^{W(t)} dt}}
    &= C_W \int_{D} \int_\R \frac{\sup_{t\in\R} e^{f(t + s) - f(s)}}{\int_\R e^{f(t + s) - f(s)} dt}  e^{f(s)} ds \mathbb P_{F_W}(df)\\
    &= C_W \int_{D} \int_\R \frac{e^{f(s)}}{\int_\R e^{f(t + s) } dt}   ds \mathbb P_{F_W}(df)\\
    & = C_W\in (0,\IF),
  \end{align*}
  where the second last equality follows from the assumption that
  $\sup_{t\in\R} F_W(t) = 0$ a.s.
  In the case $\delta > 0$ we can use the same arguments together with \netheo{prop_discr}. \QED

\section{Appendix}

The notion of a mixed moving maxima process on $\R$ defined in
\eqref{M3} can be extended to the lattice $\delta \Z$; see for instance
Remark 7 in \cite{KabExt}. Suppose that $\{\xi^\delta_ {W}(t), t \in \delta \Z\}$
is a stationary max-stable process (with standard Gumbel margins)
given by the construction \eqref{eA} with a process $W$, restricted to $\delta \Z$. Further suppose that
$W(0)=0$ almost surely and let $F_i^\delta$ be independent copies of a process $F_W^\delta$ on $\delta \Z$
with
$$\sup_{t\in \delta \Z}  F_W^\delta(t) = F_W^\delta(0) = 0$$
 almost surely, and
\begin{align}\label{Cdelta}
  C^\delta_W = \left(\E{ \sum_{t\in \delta\Z} \exp( F_W^\delta(t))}\right)^{-1} \in (0,\IF).
\end{align}
We say that $\xi^\delta_{W} $ admits an M3 representation, if
\begin{align*}
   \xiW ^\delta(t) = \max_{i\geq 0}( F_i^\delta(t - P_i^\delta) + Q_i^\delta), \quad t\in \delta \Z,
\end{align*}
where $\sum_{i=1}^\infty \epsilon_{(P^\delta_i, Q^\delta_i)}$ is a Poisson point process with intensity
$ C^\delta_W  \nu_\delta(dt) \,  e^{-x}dx$. Here $\nu_\delta / \delta$ is the
counting measure on $\delta \Z$. Below we present the counterpart of Theorem 4.1 in \cite{engelke14}
for M3 processes on lattices. We omit its proof since it follows with the same arguments as the aforementioned one.

\BT\label{prop_discr}
\EEH{Suppose that the max-stable and stationary process $\xiW$ has c\`adl\`ag sample paths.}   The process $W^\delta, \delta >0$, the restriction of $W$ to $\delta \Z$,
  can be expressed in terms of the spectral function $F^\delta$ as
  \begin{align*}
    \P(W^\delta \in L) =  C^\delta_W \mathbb{E}\left\{  \sum_{t\in\delta\Z}  \mathbf{1} \left\{ F_W^\delta(\cdot + t) - F^\delta_W(t) \in L\right\} \exp({ F^\delta_W(t)})\right\},
  \end{align*}
  which is well-defined probability measure by \eqref{Cdelta}.
\ET

\section*{Acknowledgment}
\KD{ We kindly acknowledge financial support by the Swiss National Science Foundation and
partial support by the project RARE-318984 (an FP7 Marie Curie IRSES Fellowship).
KD also acknowledges partial support 
by NCN Grant No 2015/17/B/ST1/01102 (2016-2018).}

\newcommand{\equaldis}{\stackrel{d}{=}}

\bibliographystyle{plainnat}
\bibliography{GPKF}

\end{document}